\documentclass[12pt,oneside,reqno]{amsart}
\usepackage{caption,subcaption}
\usepackage{dcolumn}
\usepackage{amsmath,amsthm,amssymb}
\usepackage{todonotes}
\usepackage{url}
\usepackage{placeins}
\usepackage{lmodern}

\usepackage{array,booktabs}

\def\be{\begin{equation}}
\def\ee{\end{equation}}

\newtheorem*{theorem*}{Theorem}
\newtheorem*{conj*}{Conjecture}
\newtheorem{conjec}{Conjecture}
\theoremstyle{remark}
\newtheorem*{rem*}{Remark}

\def\ra{\rightarrow}

\def\w{\mathsf{w}}

\begin{document}

\title[Characteristic approach to the soliton resolution]{Characteristic approach\\ to the soliton resolution}

\author{Piotr Bizo\'n}
\address{Institute of Theoretical Physics, Jagiellonian
University, Krak\'ow}
\email{piotr.bizon@uj.edu.pl}
\author{Bradley Cownden}
\address{Institute of Theoretical Physics, Jagiellonian
University, Krak\'ow}
\email{bradley.cownden@uj.edu.pl}
\author{Maciej Maliborski}
\address{University of Vienna, Faculty of Mathematics, Oskar-Morgenstern-Platz 1, 1090 Vienna, Austria, and University of Vienna, Gravitational Physics, Boltzmanngasse 5, 1090 Vienna, Austria}
  \email{maciej.maliborski@univie.ac.at}

\date{\today}%
\begin{abstract}
As a toy model for understanding the soliton resolution phenomenon we consider a characteristic initial boundary value problem for the 4$d$ equivariant Yang-Mills equation outside a ball. Our main objective is to illustrate the advantages of employing  outgoing null (or asymptotically null)  foliations in analyzing the  relaxation processes due to the dispersal of energy by radiation.  In particular, within this approach it is evident that the endstate of evolution must be non-radiative (meaning vanishing flux of energy at future null infinity). In our toy model such non-radiative configurations are given by a static solution (called the half-kink) plus an alternating chain of $N$ decoupled kinks and antikinks. We show numerically  that the configurations  $N=0$ (static half-kink) and $N=1$ (superposition of the static half-kink and the antikink which recedes to infinity) appear as generic attractors and we determine a codimension-one borderline between their basins of attraction. The rates of convergence to these attractors are analyzed in detail.
\end{abstract}
\maketitle

\section{Introduction}
According to the soliton resolution conjecture,  global-in-time generic solutions of nonlinear dispersive wave equations resolve for $t\ra\infty$  into a superposition of decoupled nonlinear bound states (solitons) and radiation \cite{T}. There are numerous physical manifestations of this phenomenon, ranging from the formation of solitons in optical fibers (modelled by a nonlinear Schr\"odinger equation) to the formation of stationary black holes in  binary black hole mergers (modelled by the Einstein equation).

 The past decade has seen a significant progress in mathematical understanding of the soliton resolution, especially for radial solutions of the energy critical nonlinear wave equations (see \cite{K} for a survey and references therein). Notably, the soliton resolution was recently proved for the
equivariant wave maps $\mathbb{R}^{2+1}\ra \mathbb{S}^2$ \cite{DKMM, JL} and the equivariant Yang-Mills (YM) equation in $4+1$ dimensions \cite{JL} (both for the global and blowup solutions).
These remarkable results are abstract in the sense that they enumerate all possible asymptotic scenarios
but do not settle which scenarios are actually realized.

As a toy model for more quantitative description of the soliton resolution phenomenon, in this paper we consider the  equivariant  YM equation in $4+1$ dimensions
\begin{equation}\label{eqr}
  W_{tt} = W_{rr} + \frac{1}{r}  W_r +\frac{2}{r^2} \, W(1-W^2),
\end{equation}
where $W(t,r)$ is the YM potential.
As the spatial domain we take the exterior of the unit ball, i.e. $r\geq 1$,  and impose the Dirichlet  condition on the boundary  $W(t,r=1)=0$.
 The associated conserved energy  is
\begin{equation}\label{E}
  E[W] = \frac{1}{2} \int_1^{\infty} \left(W_t^2 +W_r^2 + \frac{(1-W^2)^2}{r^2}  \right) r dr.
\end{equation}
Finiteness of  energy requires that $|W(t,\infty)|=1$. Since the singular point $r=0$ is outside the  domain, it is easy to see that solutions starting from smooth, finite-energy initial data ($W(0,r),W_t(0,r)$), which are compatible with the boundary condition, remain smooth for all times. Our goal is to describe their asymptotic behavior for $t\rightarrow \infty$.

 The key role in our analysis will be played by the half-kink
\begin{equation}\label{Q}
  Q(r)=\frac{r^2-1}{r^2+1}\,,
\end{equation}
which is the unique (modulo sign) static solution of equation \eqref{eqr} satisfying the boundary condition $Q(1)=0$. The half-kink is a global minimizer of energy
(see the Bogomolnyi inequality \eqref{bogom} below) and thereby a natural candidate for an attractor.
Indeed, we will see that on any compact interval $[1,R)$ every smooth solution $W(t,r)$ converges to $Q(r)$ or $-Q(r)$ as $t\ra \infty$. However, for sufficiently large energies a nontrivial coherent structure can simultaneously develop in the asymptotic region $(R,\infty)$. This behavior is intimately related to the energy criticality of the model and is absent in   supercritical dimensions; cf. the soliton resolution for the equivariant wave maps exterior to a ball in $3+1$ dimensions  \cite{BCM, KLS}.
\vskip 0.2cm
The paper is organized as follows. In section~2 we first recall from \cite{JL} the formulation of the soliton resolution conjecture for  equation \eqref{eqr} in the whole space. Then we present an analogous conjecture in our model and support one special case by  a heuristic argument based on the method of collective coordinates. In section~3 we reformulate the initial-boundary problem in terms of null foliations of constant retarded time and compactified  spatial domain. Using this formulation, in section~4 we consider the late-time behavior (the quasinormal ringdown and the polynomial tail)  for the linearized problem.
Finally, in section~5 we present  numerical evidence supporting the soliton resolution conjecture.

\section{Soliton resolution}  Let us first recall what is known about  equation \eqref{eqr} posed on the whole space $r\geq 0$. In this case, equation \eqref{eqr} is invariant under scaling, i.e. if $W(t,r)$ is a solution, so is $W_{\lambda}(t,r)=W(t/\lambda,r/\lambda)$ for any positive number~$\lambda$. The conserved energy
\begin{equation}\label{E0}
  E_0[W] = \frac{1}{2} \int_0^{\infty} \left(W_t^2 +W_r^2 + \frac{(1-W^2)^2}{r^2}  \right) r dr
\end{equation}
is scale invariant, i.e.  $E_0[W_{\lambda}]=E_0[W]$,  which is a distinctive feature of the critical dimension $d=4$ allowing for the existence of nontrivial static solutions in the presence of scaling symmetry. These static solutions, hereafter called kinks (also referred to in the literature as instantons, solitons, or bubbles), form a one-parameter family
\begin{equation}\label{Qlambda}
  Q_{\lambda}(r)=\frac{r^2-\lambda^2}{r^2+\lambda^2}
\end{equation}
with  energy $E_0[Q_{\lambda}]=\frac{4}{3}$ which is the minimum  energy  for solutions interpolating between different vacuum states $W=\pm 1$  at the origin and at infinity. Obviously, the antikink $-Q_{\lambda}(r)$ is also the solution with the same energy.
\vskip 0.2cm
 Jendrej and Lawrie proved (see Theorem~1 in \cite{JL}) that any finite-energy solution of equation \eqref{eqr} posed on the whole space tends (modulo sign) either to the vacuum  $W=1$ or to an alternating  chain of $N$  rescaled kinks and antikinks\footnote{Strictly speaking \cite{JL} deals with 2$d$ equivariant wave maps which split into equivariance classes indexed by a positive integer $k$. The case $k=2$ is essentially   equivalent to the 4$d$ equivariant YM. The only qualitative difference is that for wave maps  there are infinitely many  topological sectors, while for YM there are only two sectors (modulo sign). For this reason, in the case of wave maps   the chain of kinks and antikinks in Theorem~1 in \cite{JL} need not be  alternating.}
\begin{equation}\label{chain}
 1+\sum_{j=1}^N (-1)^{N+j} \left(Q_{\lambda_j(t)} (r) -1\right)\,.
\end{equation}
Here $\lambda_j(t)$ are continuous positive functions such that for each $j=1,...,N$
$$
\frac{\lambda_j(t)}{\lambda_{j+1}(t)} \ra 0\quad\text{as}\quad
\begin{cases}
t\ra \infty\quad \text{(for global-in-time  solutions)}\\
t\ra T \quad \text{(for blowup at finite time $T$)},
\end{cases}
$$
where by convention $\lambda_{N+1}(t)=t$ (in the  global case) or $\lambda_{N+1}(t)=T-t$ (in the blowup case), corresponding to the non-existing self-similar expansion and collapse.
\vskip 0.2cm
We return now to our toy model and make the soliton resolution conjecture:

\begin{conjec}
Any smooth, finite-energy solution $W(t,r)$  of equation \eqref{eqr} subject to the boundary condition $W(t,1)=0$ tends for $t\rightarrow \infty$  (modulo sign) either to the half-kink or to the rescaled half-kink  plus an alternating  chain of $N$ rescaled kinks and antikinks:
\begin{equation}\label{chain2}
\begin{cases}
\quad Q(r)  & \text{if $N=0$,}\\
-Q_{\mu(t)}(r)+\sum\limits_{j=1}^N (-1)^{j+1} Q_{\lambda_j(t)}(r)+1 & \text{if $N$ is odd,}\\
\quad Q_{\mu(t)}(r) + \sum\limits_{j=1}^N (-1)^{j} Q_{\lambda_j(t)}(r) & \text{if $N\geq 2$ is even.}
 \end{cases}
\end{equation}
Here $\lambda_j(t)$ are continuous positive  functions such that for each $j=1,...,N$
\begin{equation}\label{lambda}
\lambda_j(t) \ra \infty\quad \text{and}\quad  \frac{\lambda_j(t)}{\lambda_{j+1}(t)}  \ra 0\quad  \text{as} \quad t\ra\infty,
\end{equation}
where by convention $\lambda_{N+1}(t)=t$.
The function $\mu(t)$ is determined by the functions $\lambda_j(t)$ through the boundary condition $W(t,1)=0$ which implies   that $\mu(t) \ra 1$ as $t\ra \infty$.
\end{conjec}
In the  rest of the paper we confirm this conjecture for $N=0$ and $N=1$ and determine the rate of convergence to the attractors. In addition, we find a borderline between the basins of attraction using bisection along an interpolating one-parameter family of initial data.
\vskip 0.2cm
Before presenting the results of numerical simulations, we wish to put forward a heuristic argument based on the method of collective coordinates that helps to understand some aspects of asymptotic dynamics. According to \eqref{chain2}, in the case $N=1$ the attractor has the following form:
\begin{equation}\label{ansatz}
  W(t,r)=1-Q_{\mu(t)}(r)+Q_{\lambda(t)}(r), \qquad \mu^2(t)=\frac{\lambda^2(t)-1}{\lambda^2(t)+3},
\end{equation}
where the formula for $\mu(t)$ follows from the boundary condition $W(t,1)=0$.
Inserting this ansatz  into the lagrangian
\begin{equation}\label{action}
  L = \frac{1}{2} \int_1^{\infty} \left(W_t^2 -W_r^2 - \frac{(1-W^2)^2}{r^2}  \right) r dr
\end{equation}
and integrating over $r$,  in the limit of large $\lambda$ we get the effective lagrangian (we retain only  the first two leading  terms)
\begin{equation}\label{leff}
  L_{\text{eff}}=\left(\frac{4}{3}-\frac{32}{\lambda^4}\right)\, {\dot \lambda}^2 - \left(2-\frac{16}{\lambda^2}\right)\,,
\end{equation}
hence
\begin{equation}\label{Eeff}
  \left(\frac{4}{3}-\frac{32}{\lambda^4}\right)\, {\dot \lambda}^2 + \left(2-\frac{16}{\lambda^2}\right) = E_{\text{eff}}=\text{const}.
\end{equation}
Thus, the $\lambda$-particle starting at some large $\lambda(0)$ with velocity $\dot{\lambda}(0)>0$ escapes to infinity if $E_{\text{eff}}\geq 2$, while if $E_{\text{eff}}<~2$ it reaches a turning point in a finite time. The ODE \eqref{Eeff} provides a qualitative picture of the attractive interaction between the anti-half-kink and  the expanding outer kink. This approximation ceases to work when the outer kink starts shrinking because the PDE solution is no longer close to the ansatz \eqref{ansatz} (see Fig.~4 below). At the quantitative level the predictions of the effective model should be taken with caution because the ansatz \eqref{ansatz} neglects radiation. In particular, according to Conjecture~1 the expansion rate $\dot{\lambda}(t)$ in \eqref{chain2} must go to zero as $t\ra\infty$, whereas the ODE yields  $\dot{\lambda}(\infty) >0$ if $E_{\text{eff}}>2$.

\section{Characteristic formulation} We now reformulate our problem as the characteristic initial boundary value  problem. To this end we define new coordinates
\[
u=t-r, \qquad x=r^{-\frac{1}{2}}.\]
Then, on the interval $0<x\leq 1$, the YM potential $\w(u,x)=W(t,r)$ satisfies
\begin{align}\label{eq}
-4 x\, \w_{xu} + 4 \w_u  = x^4\, \w_{xx} + x^3\, \w_x +8 x^2  \w (1-\w^2),\\
\w(u,1)=0,\qquad \w(0,x)=g(x), \label{bc}
\end{align}
where the function $g(x)$ is assumed to be smooth and satisfying the finite energy condition $g(0)=1$. For such data, the solution $\w(u,x)$ remains smooth for all future times $u>0$. Moreover, the results by Chru\'sciel and collaborators \cite{CL, CW} imply\footnote{We are grateful to Piotr Chru\'sciel and Roger Tagn\'e Wafo for checking that the hypotheses of theorems on propagation of polyhomogeneity in \cite{CL, CW} hold for our equation. This, together with the absence of $\log{x}$ terms in the formal  polyhomogeneous series, shows that the solution is smooth in $x$. } that the following asymptotic expansion holds near $x=0$:
\begin{equation}\label{scri}
  \w(u,x) = 1 + \sum_{n=1} c_n(u) x^n\,.
\end{equation}
Inserting this expansion into equation \eqref{eq} and equating the coefficients of the same powers of $x$, we obtain an infinite system of ordinary differential equations for the coefficients $c_n(u)$. This system  can be solved recursively one-by-one starting from the radiation coefficient $c_1(u)$ which is free. For $n=2$ we get  $\dot c_2(u)=0$, hence  the coefficient $c_2$ is constant (so called Newman-Penrose constant). For large $n$  the nonzero coefficients $c_n(u)$ grow polynomially for $u\ra\infty$ which is a reflection of the well-known fact that the expansion \eqref{scri} is not uniform; see e.g. \cite{BF}.
\vskip 0.2cm
Multiplying equation \eqref{eq} by $x^{-3} \w_u$ we get the local conservation law
\begin{equation}\label{claw}
\partial_u \left( \frac{x}{2}  \w_x^2 + \frac{2}{x}(1-\w^2)^2\right)
     =\partial_x \left(\frac{2}{x^2}\,\w_u^2 +  x\,\w_u \,\w_x \right).
\end{equation}
Integrating this over $x$ and using \eqref{scri}, we obtain the energy loss formula
\begin{equation}\label{loss}
  \frac{d\mathcal{E}}{du} = - \dot c_1^2(u)\,,
\end{equation}
where
\begin{equation}\label{bondi}
  \mathcal{E}[\w]:=\int_{0}^{1} \left(\frac{1}{4}\,\w_x^2 + \frac{1}{x^2} ( 1-\w^2)^2\right)\,x\, dx
\end{equation}
is the Bondi-type energy (hereafter just called energy).
Note that $\mathcal{E}[\w]$ is equal to the potential part of the total conserved energy $E[W]$.
In terms of $x$ the half-kink reads
\begin{equation}\label{q}
  q(x):=Q\left(\frac{1}{x^2}\right)=\frac{1-x^4}{1+x^4}.
\end{equation}
It is the  global minimizer of $\mathcal{E}[\w]$ as follows from the Bogomolnyi inequality
 \begin{align}\label{bogom}
 \mathcal{E}[\mathsf{w}] = \int_0^1 \left[\frac{1}{2}\w_x+\frac{1}{x} (1-\w^2)\right]^2\, x\, dx - \int_0^1 \partial_x (\w-\frac{1}{3} \w^3)\,dx \geq  \frac{2}{3},
 \end{align}
 which is saturated on $\w=q(x)$, i.e. $\mathcal{E}[q]=\frac{2}{3}$.
\vskip 0.2cm
Since  $\mathcal{E}[\w]$ is non-increasing and bounded below, there exists a limit $$\mathcal{E}_{\infty}=\lim_{u\rightarrow \infty} \mathcal{E}(\w(u))\geq \mathcal{E}(q)=\frac{2}{3}\,.$$
According to Conjecture~1 the endstate (which clearly  must be non-radiative) has the form \eqref{chain2}, hence the final energy is a sum of the energies of the half-kink and $N$ kinks/antikinks
\begin{equation}\label{EN}
 \mathcal{E}_{\infty}=\frac{2}{3}+\frac{4}{3} N.
 \end{equation}
In section~5 we describe the relaxation to  the $N=0$ and $N=1$ attractors. In our numerical simulations we have not observed $N\geq 2$ attractors which suggests that, if they exist,  they are nongeneric.

\section{Linearized dynamics near the half-kink}
Let $\w=q(x)+x f(u,x)$. Substituting this into equation \eqref{eq} we obtain
\begin{equation}\label{eqf}
   f_{ux} + \frac{1}{4}\partial_x (x^3 f_x) - U(x) f - 6 q(x) x^2 f^2 - 2 x^3 f^3 =0,\qquad f(u,1)=0,
\end{equation}
where
\begin{equation}\label{V}
  U(x)=\frac{x (15-66 x^4 +5 x^8)}{4(1+x^4)^2} = \frac{15}{4} x +\mathcal{O}(x^5).
\end{equation}
Dropping the nonlinear terms and the $\mathcal{O}(x^5)$ term in the potential, we get the linear equation (corresponding to linearization around $\w=1$  rather than $q$)
\begin{equation}\label{eqfl}
  f_{ux} + \frac{1}{4} \partial_x (x^3 f_x) - \frac{15}{4} x f=0.
\end{equation}
This equation has an explicit solution
\begin{equation}\label{sol_f}
  f_0(u,x)= u^{-\frac{5}{2}} (u x^2 +2)^{-\frac{5}{2}}.
\end{equation}
General solutions of equation \eqref{eqfl} for initial data with  the vanishing NP constant  behave similarly  to $f_0(u,x)$ for $u\ra \infty$, i.e. they decay as $u^{-5}$ in the interior  ($x>0$) and  $u^{-5/2}$ at future null infinity ($x=0$). This can be shown directly, or  by defining $F(t,r)=x^5 f(u,x)$ and rewriting equation \eqref{eqfl} in terms of the original coordinates $(t,r)$. Then, $F(t,r)$ satisfies  the free radial wave equation in $6+1$ dimensions $F_{tt} - F_{rr} - \frac{5}{r}  F_r=0 $, for which the late-time pointwise decay $F(t,r)\sim (t-r)^{-5/2} (t+r)^{-5/2}$ is well known \cite{H}. The rate of decay of  linear perturbations about $q(x)$ is the same because the term $\mathcal{O}(x^5)$ in the potential \eqref{V} is asymptotically negligible. The numerical verification of this claim is depicted in Fig.~1 where we also
show that solutions with nonzero NP constant exhibit slower decay.
\begin{figure}[h]
	\includegraphics[width=0.50\textwidth]{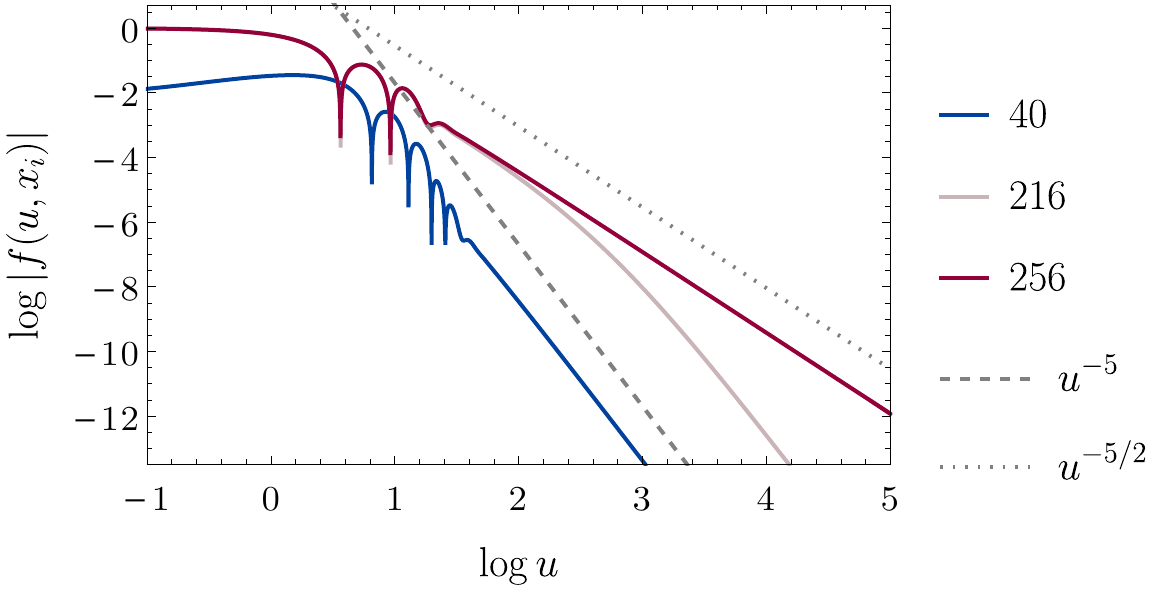}
	\includegraphics[width=0.47\textwidth]{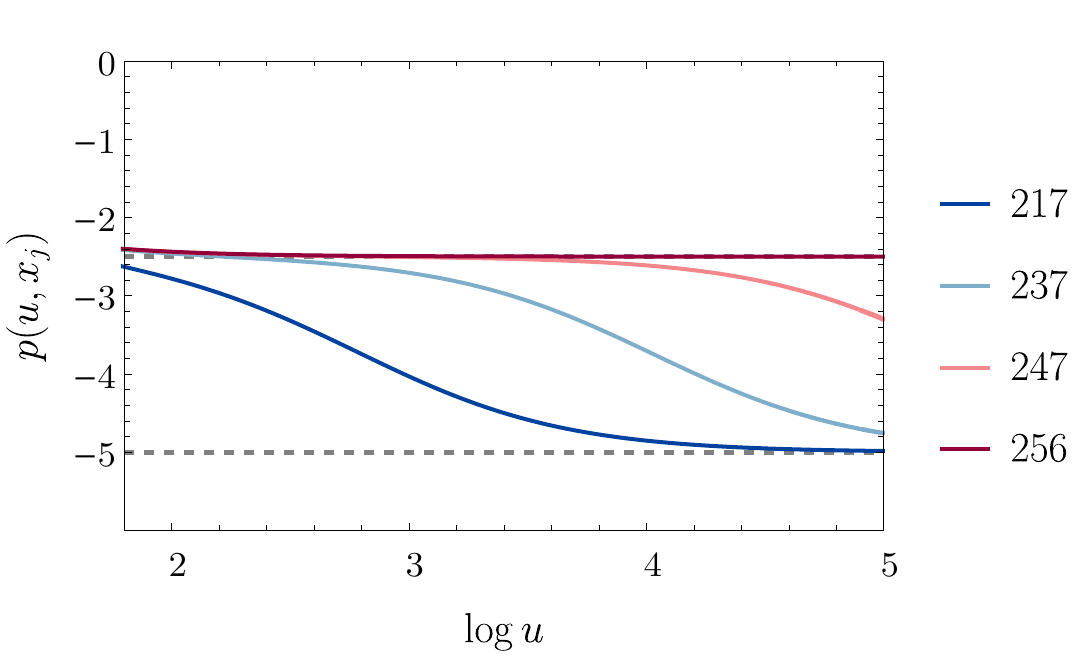}
	\vspace{2ex}
	\includegraphics[width=0.50\textwidth]{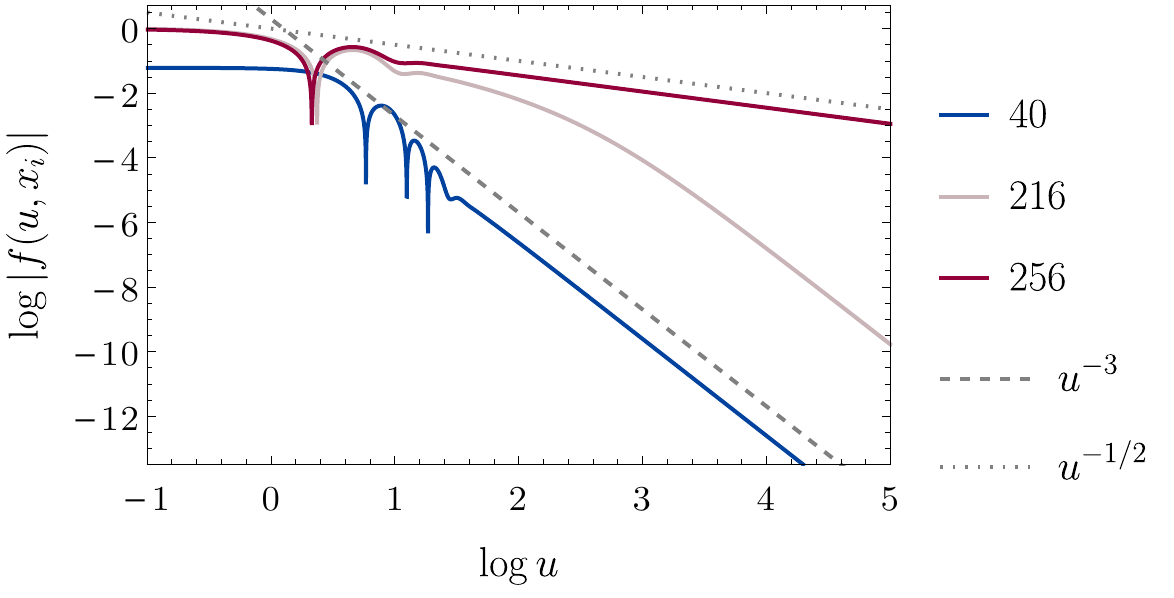}
	\includegraphics[width=0.47\textwidth]{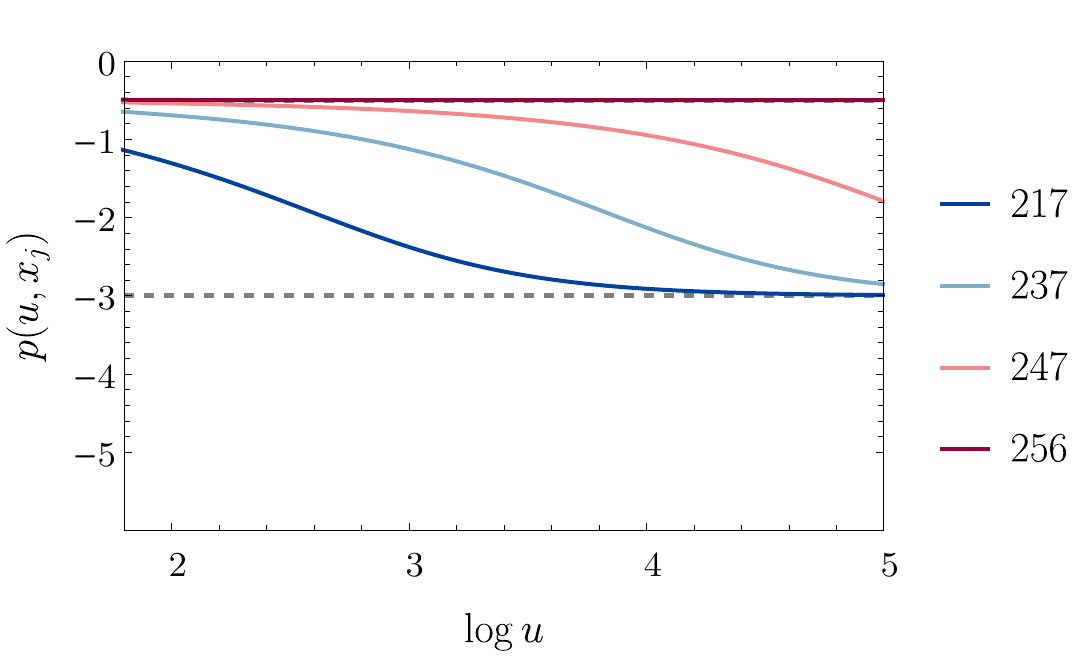}
	\caption{\small{Linear evolution about the half-kink $q(x)$ for sample initial data with zero and nonzero NP constant, respectively: $f(0,x)=\cos^{2}\!\left(\pi x/2\right)$ (top row) and  $f(0,x)=\cos^{2}\!\left(\pi x/2\right)+(1-x)x$ (bottom row). On the left, we plot $\log|f(u,x_i)|$ at different grid points (solid lines labelled by the grid point number $i$, where $x_{256}=0$; see section~5 for the details) together with the theoretical  decay rates (dashed/dotted lines). On the right, we plot the local power index $p(u,x_j):=u\,\partial_{u}f(u,x_i)/f(u,x_i)$ evaluated at different grid points.}
	}
	\label{fig:linear_pointwise_decay_q}
\end{figure}
\newpage
 Having determined the late-time linear tail, now we turn to the computation of quasinormal modes. They govern the relaxation to the half-kink for intermediate times before the tail is uncovered.  It is convenient to rewrite the linear part of equation \eqref{eqf} in terms of $y=x^2$. Substituting
 \begin{equation}\label{vf}
 f(u,x)=e^{s u} \,v(y),
 \end{equation}
 we get the eigenvalue problem
 \begin{equation}\label{eqv}
 2s v' +(y^2 v')' - V(y) v =0, \qquad V(y)=\frac{15}{4}-\frac{48+5y^2}{(1+y^2)^2}\,y^2
\end{equation}
with $v(1)=0$. Following Leaver \cite{L} we seek solutions  of \eqref{eqv} in terms of the power series
\begin{equation}\label{power}
  v(y)=\sum_{n\geq 1} a_n (1-y)^n,\qquad a_1=1.
\end{equation}
Since the nearest singularity from $y=1$ is located at $y=0$, this power series is absolutely convergent for $y\in (0,1]$. The eigenvalues $s_n$, called quasinormal frequencies, are selected by the condition that the power series is absolutely convergent at $y=0$; as follows from \eqref{vf}, the corresponding  solutions $f_n(u,y)$, called the quasinormal modes,  are purely outgoing\footnote{An alternative way of selecting the outgoing solution by a certain Gevrey-class regularity condition has been recently proposed by Gajic and Warnick \cite{GW}; see also \cite{GZ}. For a very interesting discussion of hyperboloidal approach to quasinormal modes and Leaver's method we refer the reader to \cite{AM, MJA}.}. Inserting \eqref{power} into equation \eqref{eqv} we get a 7-term recurrence relation. Among its six linearly independent solutions, four solutions  decay as $a_n\sim  2^{-\frac{n}{2}}$ for $n\ra\infty$, hence they do not affect the convergence properties of the series at $y=0$. Using the method of successive approximations \cite{WL} one can show that the remaining two solutions have the following asymptotic expansions
\begin{equation}\label{asym_an}
  a^{(+)}_n \sim n^{-\frac{3}{4}} e^{\sqrt{8s n}} \sum_{k=0} c^{(+)}_k n^{-\frac{k}{2}},\qquad  a^{(-)}_n \sim n^{-\frac{3}{4}} e^{-\sqrt{8s n}}  \sum_{k=0} c^{(-)}_k n^{-\frac{k}{2}},
\end{equation}
where the coefficients $c^{(\pm)}_k$ can be determined successively by plugging the expansions \eqref{asym_an} into the recurrence relation. We conclude that for $n\ra \infty$
\begin{equation}\label{sol_an}
  a_n = C_+(s) a^{(+)}_n +C_{-}(s) a^{(-)}_n +\mathcal{O}(2^{-\frac{n}{2}}).
\end{equation}
For $|\arg(s)|<\pi$, the series $\sum |a^{(+)}_n|$ diverges while the series $\sum |a^{(-)}_n|$ converges, therefore
the quasinormal frequencies are given by the roots of the coefficient $C_+(s)$. There are several alternative ways to find these roots. The most frequently used  is Leaver's method of continued fractions \cite{L}. It is stable and accurate but tedious in the case at hand because  the recurrence relation must first be reduced to three terms by Gaussian elimination. Employing this method we found exactly one pair of complex conjugate frequencies
 $s  \approx-0.364322 \pm  0.476858 i$. To verify this result, we have reproduced it by two different methods: a brute force evaluation of a dominant solution by forward recurrence and an algebraic method introduced in \cite{bm}.  We skip the details of these straightforward but dull computations. We confirmed the above perturbative analysis by the direct numerical integration of the linearized equation; see Fig.~2.

 \begin{figure}[h]
	\includegraphics[width=0.60\textwidth]{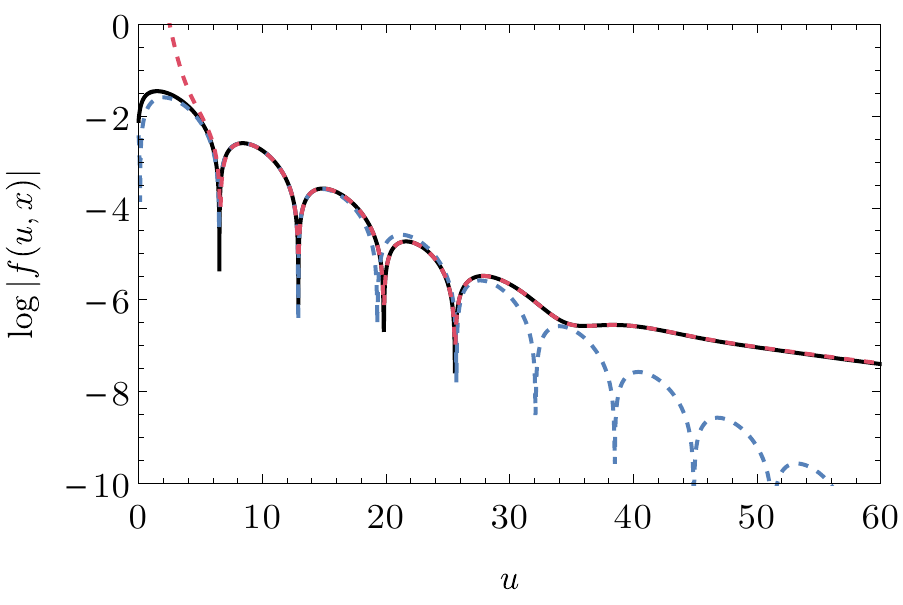}
	\caption{\small{The solution depicted in  blue on the top left plot in Fig.~1 is shown here for early times when  the relaxation to the half-kink is governed by the quasinormal mode.  The fit (dashed red line) of a superposition of an exponentially damped oscillation and the polynomial tail to the numerical data  (black line) gives the quasinormal  frequency $s=-0.364271 \pm 0.476856 \,i$, in very good agreement with the perturbative calculation. The pure quasinormal ringdown is depicted for reference by the dashed blue line.}}
	\label{fig:linear_pointwise_decay_q_qnm_fit}
\end{figure}

\begin{rem*}
\small{It is instructive to compare  the above computation of the quasinormal modes for the half-kink  with an analogous computation for the vacuum solution $\w=1$ of equation \eqref{eq} with the boundary condition $\w(u,1)=1$ (as mentioned above, this is equivalent to the free wave equation in $6+1$ dimensions). For the ansatz $\w=1+\sqrt{y}\, e^{s u} v(y)$ (where $y=x^2$), we  obtain the same eigenvalue equation as \eqref{eqv} but with the potential $V=15/4$. Repeating the above analysis, we get  a three-term recurrence relation having two linearly independent solutions $a^{(\pm)}_n$ with the asymptotic behavior \eqref{asym_an}, hence as before the quasinormal frequencies are given by the roots of the coefficient $C_+(s)$ of the dominant solution.  We remark that in this case the analysis based on the recurrence relation is purely academic  because  the eigenvalue equation can be solved exactly and the outgoing solution is given by  $v_{\text{out}}=y^{-\frac{1}{2}} e^{s/y} K_2(s/y)$, where $K_2(z)$ is the modified Bessel function of the second kind. Thus, the quantization condition for the quasinormal modes is $K_2(s)=0$, which has exactly one pair of complex conjugate zeros on the principal branch $s=-1.281373 \pm 0.4294849 i$ \cite{P}. We verified that the roots of $C_+(s)$ are the same, which provides  a reassuring benchmark test for Leaver's method.}
\end{rem*}
\section{Numerical results}
\begin{figure}[h]
	\includegraphics[width=0.50\textwidth]{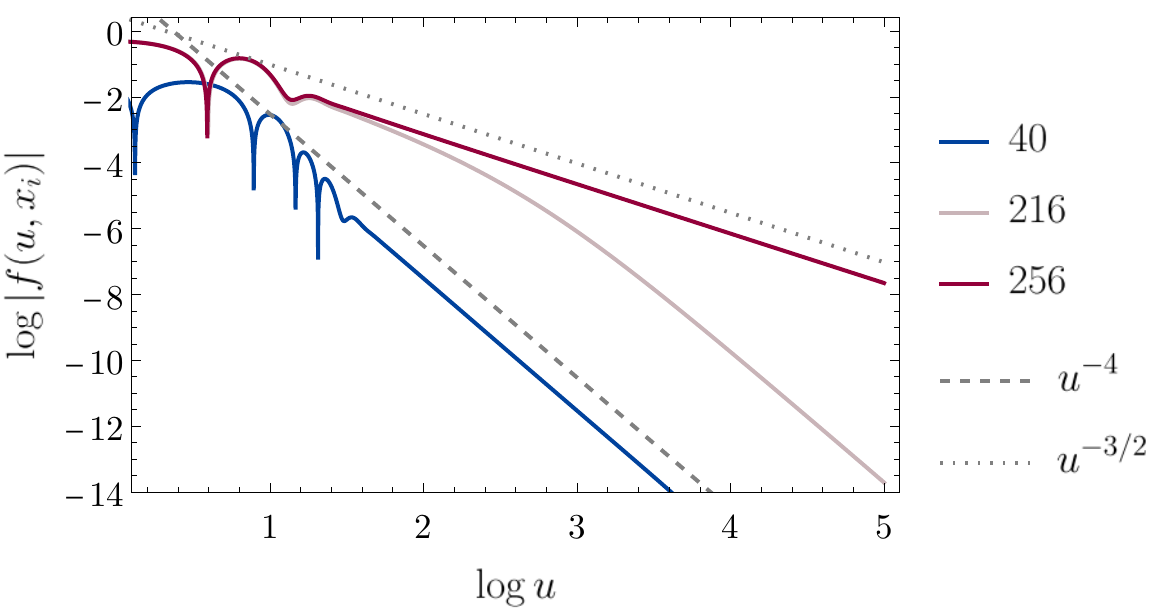}
	\includegraphics[width=0.47\textwidth]{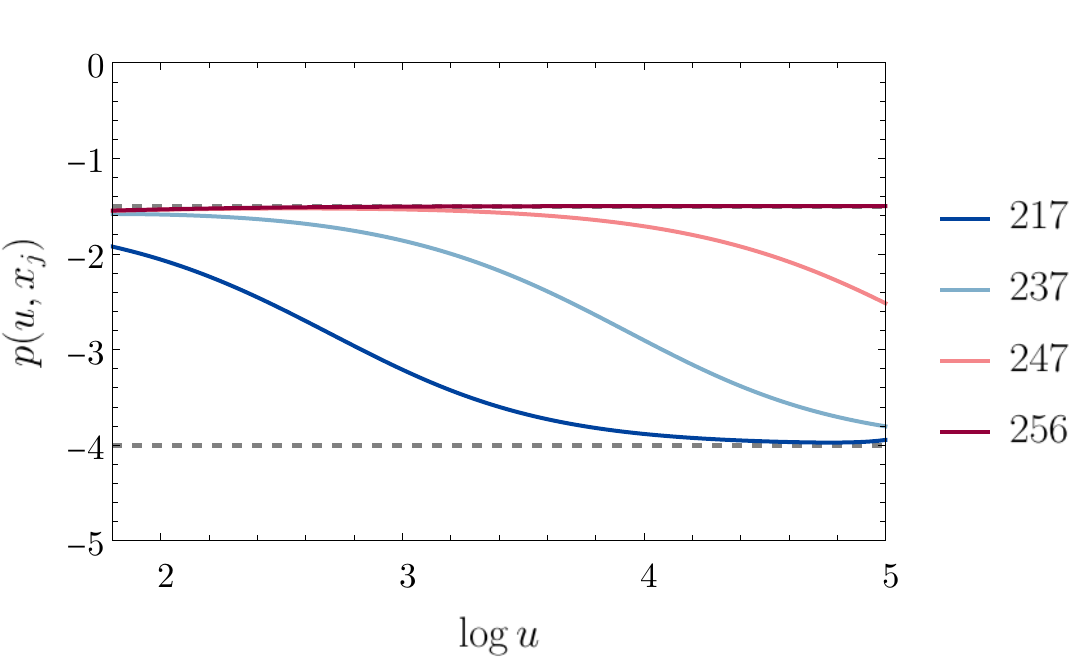}
	\vspace{2ex}
	\includegraphics[width=0.50\textwidth]{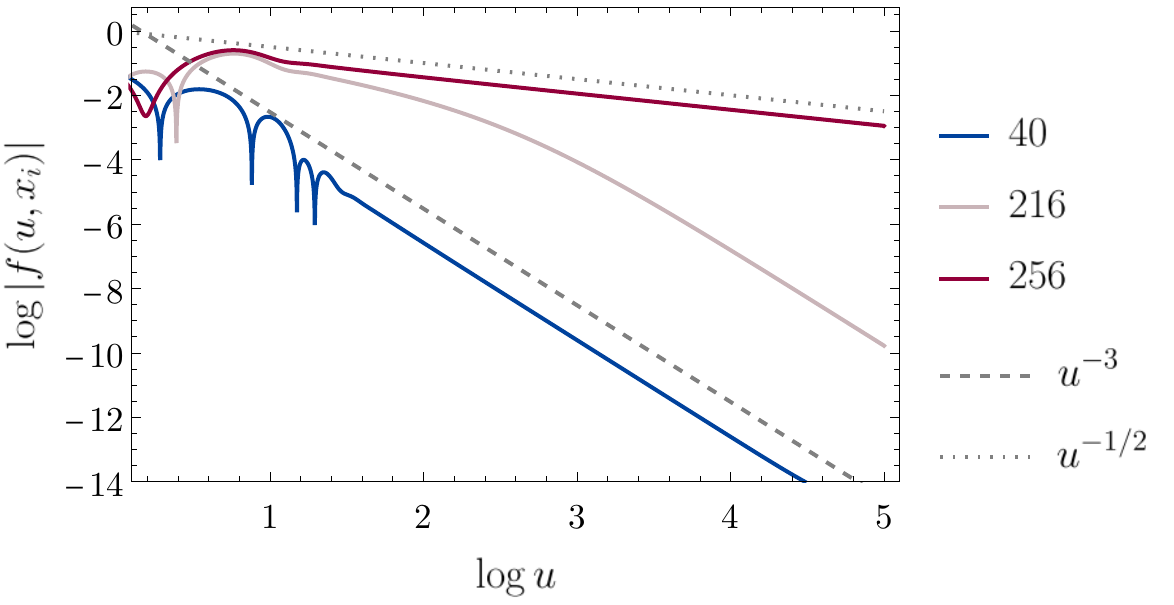}
	\includegraphics[width=0.47\textwidth]{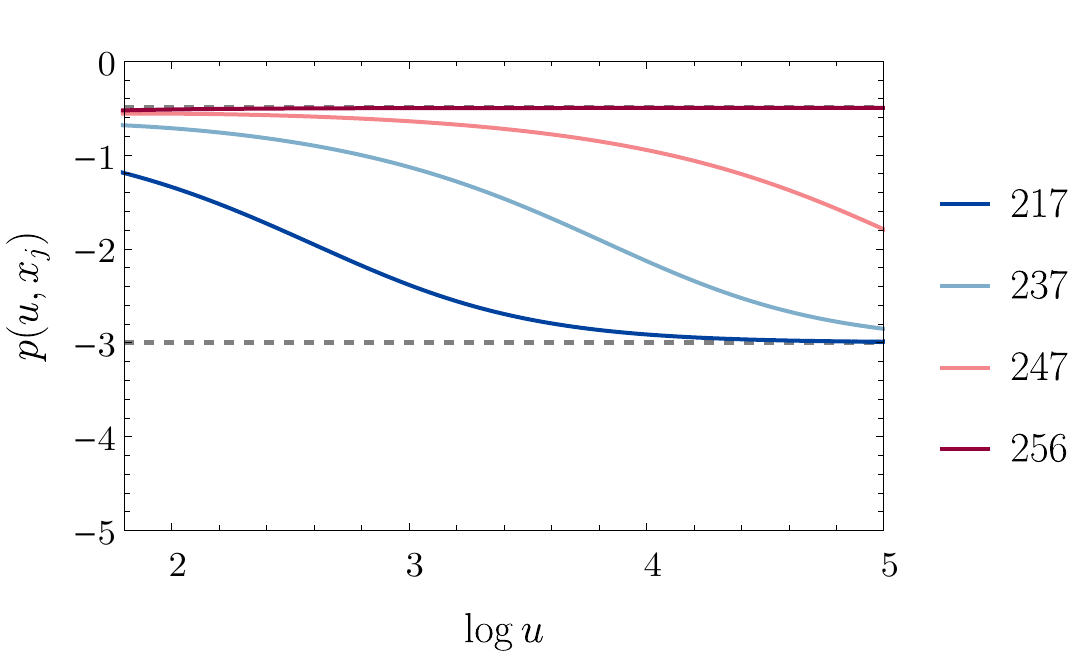}
	\caption{\small{Pointwise decay to the half-kink  for the initial data~\eqref{icb} with $b=1$ (top row) and for the same data with an extra term $x(1-~x)$ which generates a nonzero NP constant (bottom row). Compare with the analogous plots of the linear decay shown  in Fig.~\ref{fig:linear_pointwise_decay_q}.}
	}
	\label{fig:nonlinear_pointwise_decay_q}
\end{figure}

In this section we corroborate Conjecture~1 with direct numerical simulations of the initial boundary value problem \eqref{eq}-\eqref{bc}. As in section~4, we write $\w=q(x)+x f(u,x)$ and then solve equation \eqref{eqf} numerically using the method of lines. To this end, we first discretize equation \eqref{eqf} in space using the pseudospectral approach. For numerical convenience, we rescale the spatial domain to the interval $[-1,1]$ using $z=2x-1$ and
  work with function values $\{f_{j}(u)\equiv f(u,z_{j})\}$ evaluated at $K$ Chebyshev points of the second kind $\{z_{j}=\cos\left(\frac{(j-1)\pi}{K-1}\right)\}$, $1\leq j\leq K$. Spatial derivatives  $\partial_{z}$ and $\partial_{z}^{2}$ are replaced by the corresponding spectral differentiation matrices $D^{(1)}_{K}$ and $D^{(2)}_{K}$ \cite{Trefethen.2000} and then both the derivatives and nonlinear terms are evaluated using the grid function $\{f_{i}\}$. The resulting semi-discrete system takes the following schematic form
\begin{equation}
	\label{eq:28.11.21_02}
	\partial_u f_{1}=0,\qquad \partial_u(D^{(1)}_{K} f)_{j} = F_j(D^{(2)}_{K} f, D^{(1)}_{K} f, f,z)\,, \quad j=2,\ldots,K,
\end{equation}
where the boundary condition $f_1=0$ replaces the equation at the grid point  $z_{1}=1$. We bring this system to an explicit form by solving the linear system for $u$ derivatives of $\{f_j\}$. This requires inverting the operator $D^{(1)}_{N}$ with the first row replaced by a condition  $\partial_u f_{1}=0$. Note that this linear operator is invertible and non-degenerate. The solution uses the LU decomposition of the resulting matrix.
For efficient time integration we use an implicit scheme. We employ the BDF method (variable-order, variable-coefficient, in fixed-leading-coefficient form) which, for the best performance, we limit to the second order (higher-order methods struggled to find the optimal step size/order, probably due to the stiffness of the equation). We used the IDA code \cite{Hindmarsh.1999}, as available in Wolfram Mathematica \cite{Mathematica.2021}, in which we set the error tolerances to very conservative values (typically $10^{-13}$) so that the spatial resolution determines the errors in the numerical solution. Tests of the final algorithm show the spectral (exponential) convergence with increasing $K$.
\vskip 0.2cm
Using the above method we have simulated the evolution of various initial data. Here we illustrate the results for a sample one-parameter family
\begin{equation}
	\label{icb}
	\w(0,x) = 1 + b x^4 - (1+b) x^6\,,
\end{equation}
where $b$ is a free parameter. For this data the NP constant is equal to zero. To see how  a nonzero NP constant affects the dynamics we look in parallel at the evolution of initial data \eqref{icb} with the additional term $x(1-x)$. The  energy of initial data \eqref{icb} attains the minimum value $0.6672$ at $b\approx-2.1022$.
For $-7.7295\lesssim b\lesssim 2.5933$ we have $\mathcal{E}<2$, hence according to \eqref{EN} the half-kink $q$ is the only possible attractor (the $N=0$ case in our terminology). In agreement with this,
we observe rapid convergence to $q$ through a short ringdown and then the late-time tail; see Fig.~\ref{fig:nonlinear_pointwise_decay_q}. Notice that the linear decay rates determined in the previous section, namely $u^{-5}$ in the interior ($x>~0$) and $u^{-5/2}$ along the future null infinity ($x=0$), are not propagated by the nonlinear flow\footnote{Interestingly, if the  NP constant is nonzero, then the linear and nonlinear tails decay at the same rate: $u^{-3}$ for $x>0$ and $u^{-1/2}$ for $x=0$; see the bottom rows in Fig.~1 and~3.} for which the decay rates are slower by one power of~$u$. Most important for us is the decay rate for the radiation coefficient $c_1(u) \sim u^{-3/2}$, however establishing this fact rigorously is a task that goes beyond the scope of this paper. We remark that similar nonlinear tails (but only in the interior) have been studied in the literature for semilinear wave equations in high even spatial dimensions (in particular, for the quadratic wave equation in $6+1$ dimensions which is relevant in our context); see \cite{AKT} and \cite{HV}.

Next, we consider  initial data with energy greater than $2$.   For moderate positive values of $b$   the solution  again quickly converges to the half-kink, however for larger values of $b$ we observe  formation of the superposition of the anti-half-kink and the expanding kink (the $N=1$ attractor in our terminology). We find that the transition between these two scenarios occurs at $b_{0}\approx 12.458288341217909$. For marginally subcritical solutions (i.e. for $b=b_{0}-\varepsilon$ with small positive~$\varepsilon$) a superposition of the anti-half-kink and the expanding kink appears for intermediate times but at a later time  the expansion stops, the kink starts shrinking and is quickly annihilated.
In this process the energy of the kink is rapidly radiated away and the solution settles down to the half-kink. This behavior is shown in Figs.~4~and~5 where we plot the snapshots of marginally subcritical and supercritical solutions and the corresponding energies, respectively.
\vspace{-0.5cm}
\begin{figure}[h]
	\includegraphics[width=1.\textwidth]{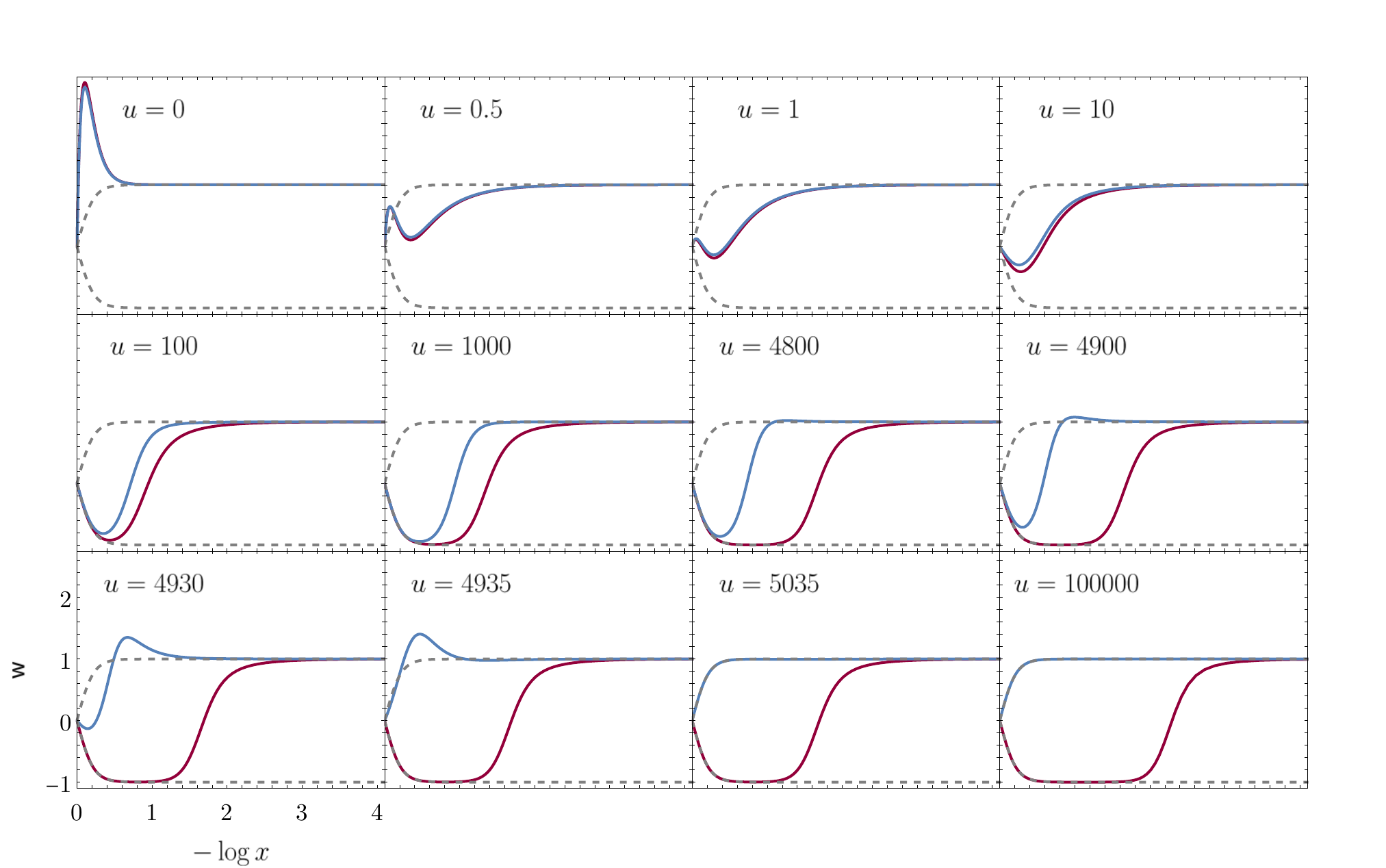}
	\caption{\small{Snapshots from the evolution for the initial data \eqref{icb} with marginally subcritical (blue) and supercritcal (red) values of $b$ near $b_{0}$. The half-kink and anti-half-kink are plotted with dashed lines. For the subcritical evolution the turning point is at $u_R \approx 1000$. }}
	\label{fig:nonlinear_snapshots}
\end{figure}
\vspace{-0.5cm}
\begin{figure}[h]
	\includegraphics[width=0.5\textwidth]{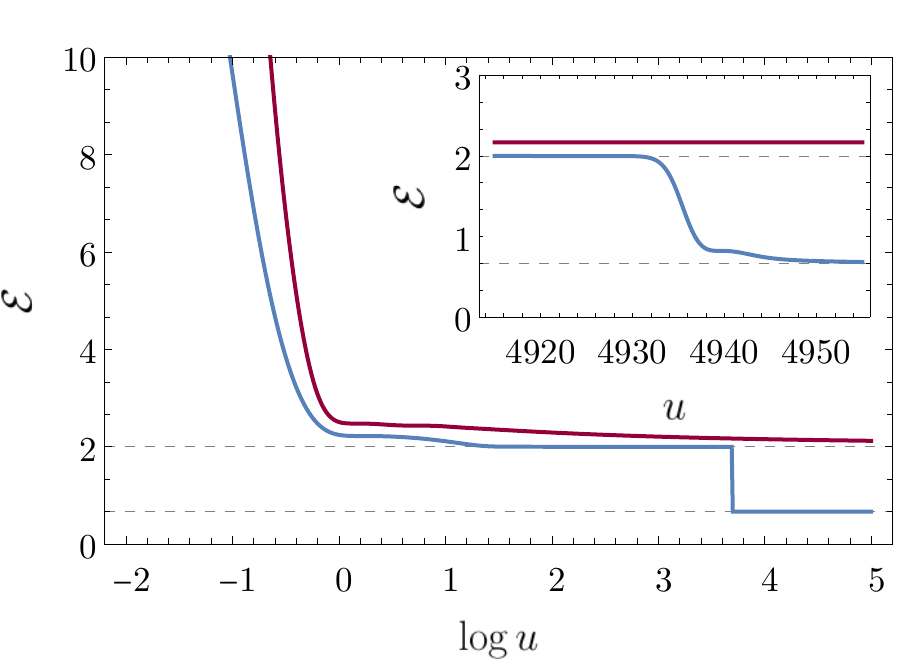}
	\caption{\small{The energies of solutions from Fig.~4. In the inset we zoom into the intermediate phase of subcritical evolution when the energy of the outer kink is radiated away.  }}
	\label{fig:nonlinear_snapshots_energy}
\end{figure}

\begin{figure}[h]
\includegraphics[width=0.49\textwidth]{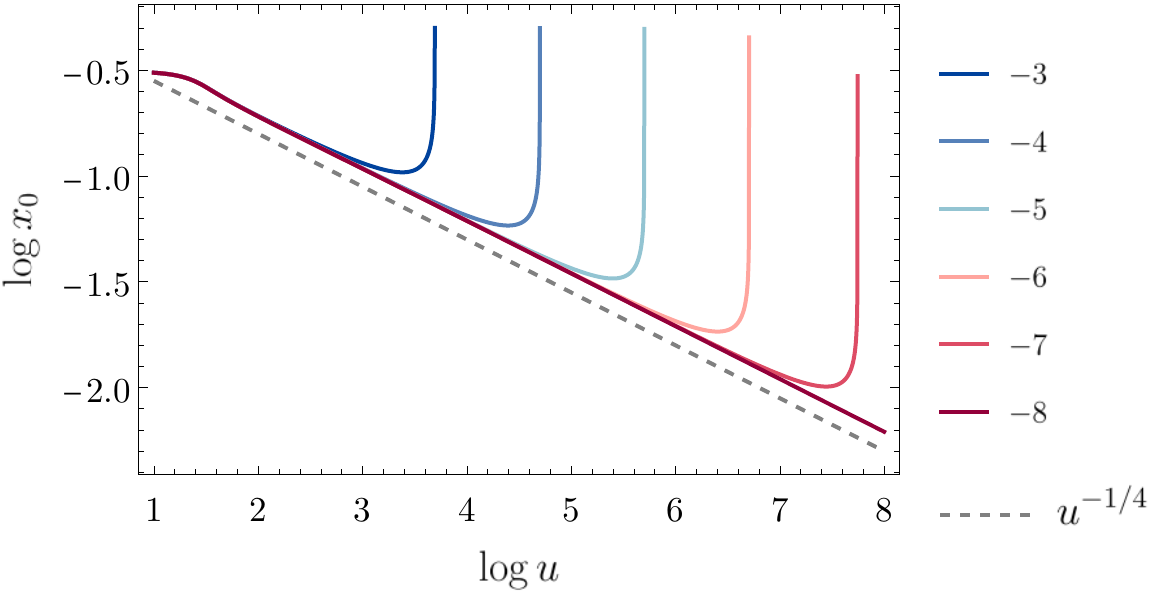}
	\includegraphics[width=0.39\textwidth]{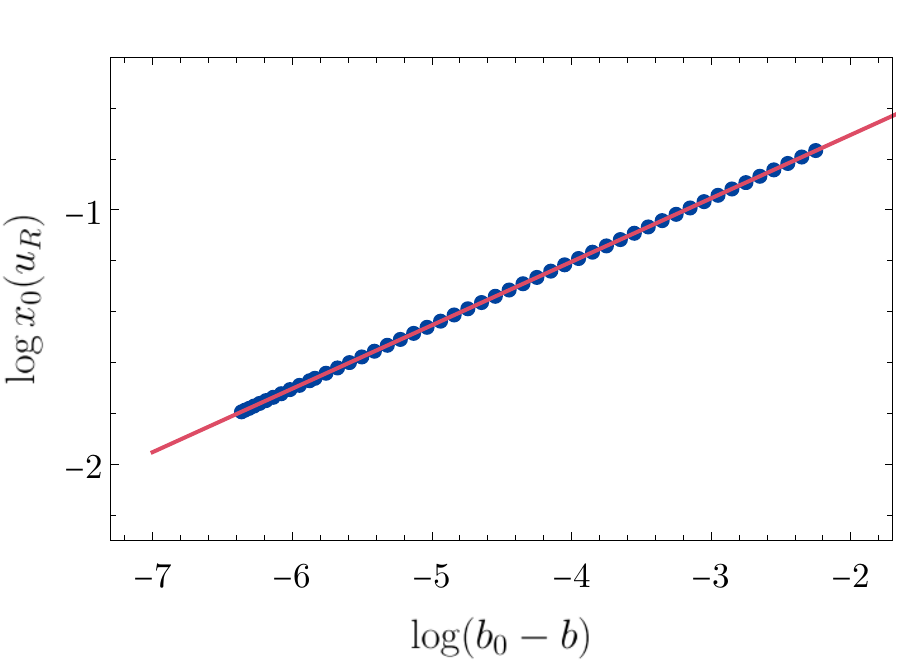}
	\caption{\small{Left panel: for marginally subcritical solutions we plot the zero of $\w(u,x)$, denoted by $x_0(u)$, for several values of $\log(b_0-b)$ (depicted by different colors). For almost critical data the fit gives $x_0\sim u^{-0.2494}$. This is in accord with the ODE approximation \eqref{Eeff} which gives $\lambda(t) \sim t^{1/2}$ for motion on the separatrix with $E_{\text{eff}}=2$. Right panel: $x_0(u)$ at the return time $u_R$ as the function of $b_0-b$. The fit gives $x_0(u_R)\sim (b_0-b)^{0.2495}$ (red line). }}
	\label{fig:subcritical_x0}
\end{figure}

For a more quantitative description of the expanding phase of subcritical solutions, let $x_0(u)$ be the zero of the solution $\w(u,x)$ and $u_R$ be the return time when the expansion stops. We find that
$u_R\sim \varepsilon^{-1}$ and $x_0(u_R) \sim \varepsilon^{1/4}$; see Fig.~\ref{fig:subcritical_x0}.
Translating these scaling relations to the original variables and comparing with the ansatz \eqref{ansatz} we get $t_R\sim \varepsilon^{-1}$ and $\lambda(t_R) \sim \varepsilon^{-1/2}$. This is in  agreement with the ODE approximation \eqref{Eeff} which predicts  that for marginally subthreshold effective energies we have $\lambda(t_R) \sim (2-E_{\text{eff}})^{-1/2}$.

\begin{figure}[h]
	\includegraphics[width=0.6\textwidth]{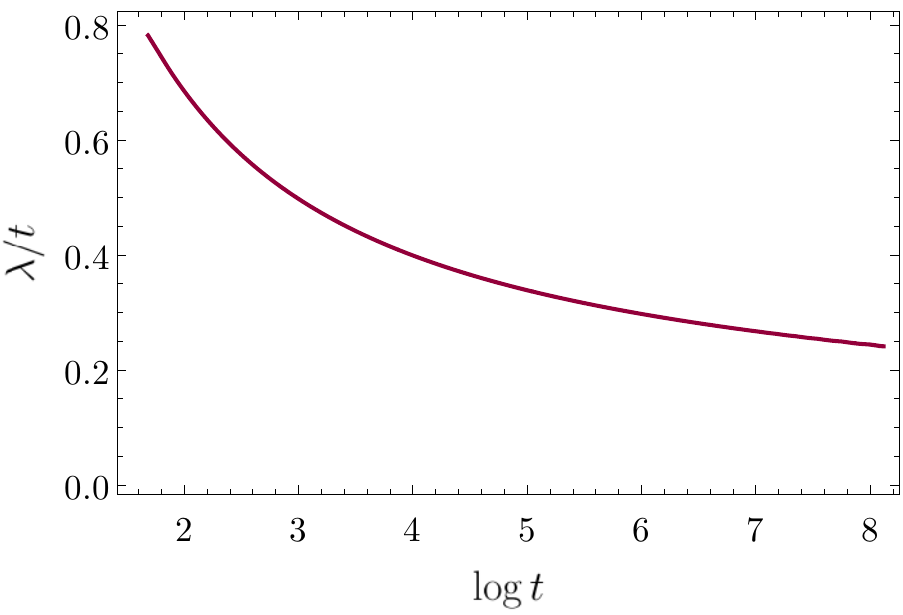}
	\caption{\small{The speed of expansion of the outer kink  for supercritical evolution with $b=20$. To compute the function $\lambda(t)$, we first find the zero $x_0(u)$ of the solution $\w(u,x)$, then translate the result to the variables $(t,r)$, and finally compare it with the ansatz \eqref{ansatz}.}}
	\label{fig:lambda_over_t_vs_t}
\end{figure}

Increasing $b$  we find that above $b_{1}\approx 47.90418049175238$ the solution again settles down on the half-kink  after a few rapid nonlinear oscillations. A similar transition occurs below $b_{-1}\approx -53.94479194728988$. As $|b|$ grows further the endstate keeps flipping back and forth between the $N=0$ and $N=1$ attractors. We conjecture that there are infinitely many critical values $b_{n}$ ($n\in\mathbb{Z}$) at which the curve of initial data \eqref{icb} intersects  the $N=0$ and $N=1$ basins of  attraction.
\vskip 0.2cm
In the case $N=1$ of Conjecture~1, it remains to verify that the speed of expansion of the outer kink goes asymptotically to zero, i.e.  $\frac{\lambda(t)}{t}\ra 0$ for $t\ra\infty$. This is shown in Fig.~\ref{fig:lambda_over_t_vs_t}.
Unfortunately, we are not in position to say more about the dynamics of $\lambda(t)$.
 The computation of the precise asymptotic behavior of $\lambda(t)$, which takes into account the loss of energy by radiation, is a challenging open problem that we leave  to future work\footnote{See \cite{BOS, RR} for the derivation of the modulation equation for $\lambda(t)$  for the blowup solutions of the $4+1$ YM equation in the whole space.}.

\vskip 0.5cm
\noindent\emph{Acknowledgement.} PB wishes to thank Piotr Chru\'sciel, Peter Hintz and Jacek Jendrej for useful comments. Special thanks are due to Arthur Wasserman for very helpful suggestions.
 The work of PB and BC was supported  by the National Science Centre grant no.\ 2017/26/A/ST2/00530.
MM acknowledges the support of the Austrian Science Fund (FWF), Project P 29517-N27 and the START-Project Y963.

\end{document}